\documentclass[11 pt,reqno]{amsart}
\usepackage{amssymb,amsmath,amstext,amsgen,amsfonts,amsthm,verbatim,hyperref,fullpage,enumitem}
\usepackage[usenames,dvipsnames]{color}
\usepackage{amsrefs}
\usepackage{tikz}
\usetikzlibrary{calc,arrows.meta}

\newcommand{\N}{\mathbb{N}}
\newcommand{\R}{\mathbb{R}}

\renewcommand{\H}{\mathbb{H}}

\renewcommand{\epsilon}{\varepsilon}
\renewcommand{\phi}{\varphi}

\newcommand{\cH}{\mathcal{H}}

\newcommand{\MCP}{\mathsf{MCP}}
\newcommand{\CD}{\mathsf{CD}}
\newcommand{\RCD}{\mathsf{RCD}}
\newcommand{\m}{\mathfrak{m}}
\DeclareMathOperator{\Geo}{Geo}

\newtheorem{theorem}{Theorem}

\theoremstyle{remark}

\begin{document}

\title{The measure contraction property and rectifiability}
\author{Sean Li}
\date{\today}
\address{Department of Mathematics, University of Connecticut, Storrs, CT 06269}
\email{sean.li@uconn.edu}

\begin{abstract}
  We show that a metric measure space $(X,d,\m)$ with finite upper $N$-density that satisfies $\MCP(K,N)$ is $N$-rectifiable. This removes a unique tangent requirement from \cite{MMR}.
\end{abstract}

\maketitle

\section{Introduction}

In the last two decades, a synthetic notion of Ricci curvature bounds has been developed around the curvature-dimension condition $\CD(K,N)$, introduced in the seminal works of Lott-Villani \cite{LV} and Sturm \cite{sturm1,sturm2}. Roughly speaking, their definition requires a $(K,N)$-convexity of an entropy functional along Wasserstein geodesics. Another condition, the measure contraction property $\MCP(K,N)$ was introduced independently (in slightly different forms) by Ohta \cite{ohta} and Sturm \cite{sturm2}. We will be concerned with Ohta's $\MCP$, which requires a $\CD$ like inequality, but only for transports from an arbitrary measurable set to a Dirac mass. In particular, $\CD(K,N)$ implies Ohta's $\MCP(K,N)$ \cite{rajala} but not the other way around as the Heisenberg group is an example of a $\MCP(0,5)$ space that is not $\CD(K,N)$ for any $K,N$ \cite{juillet}.

Much work has been devoted to the local structure of these spaces, in particular whether they are rectifiable. The Heisenberg group gives a negative example of a $\MCP$ space that is not rectifiable. On the other hand, Mondino and Naber showed that spaces satisfying a stricter condition $\RCD(K,N)$ (that we will not cover here) are rectifiable \cite{MN}. For $\MCP$, Magnabosco-Mondino-Rossi \cite{MMR} showed that $\MCP(K,N)$ spaces are rectifiable assuming finite upper $N$-density and almost everywhere unique tangents. We show that this latter condition is not needed.

Recall that $\Theta^{N*}(x) = \limsup_{r \to 0} \frac{\m(B(x,r))}{r^N}$.

\begin{theorem}
  Let $N \in (1,\infty)$ and $(X,d,\m)$ be a $\MCP(K,N)$ space and
  \begin{align}
    \Theta^{N*}(x) < \infty, \qquad \m-a.e. ~x \in X. \label{eq:density-bounds}
  \end{align}
  Then $N \in \N$ and $X$ is $N$-rectifiable.
\end{theorem}

The result in \cite{MMR} is stated with a positive lower density bound, but this is automatic for $\MCP$ spaces due to Bishop-Gromov.

Note that the Heisenberg group $\H$ shows that the dimensions between the upper density and $\MCP$ must match for such a statement. Indeed, $\dim_H \H = 4$ but, as mentioned, $\H$ satisfies $\MCP(0,5)$. The converse also fails. On $N$-dimensional Riemannian manifolds equipped with the volume measure, $\MCP(K,N)$ is equivalent to $\mathrm{Ric} \geq K$ \cite{ohta,sturm2}. Thus, any complete Riemannian $N$-manifold whose Ricci curvature is unbounded from below is an $N$-rectifiable space with finite $N$-density that satisfies no $\MCP(K,N)$.

The basic idea is to use a combination of results of \cite{BKO,BL,cheeger,EB} to say that an unrectifiable $\MCP$ space must admit Lipschitz maps $\{\phi_i\}$ into a lower dimensional Euclidean space that allow a Cheeger-style differentiability of real-valued Lipschitz functions. As the map is into a lower dimensional space, for every point there must exist many other points close to the fibers of $\phi_i$. Differentiating distance functions would then mean that points on the same fiber do not lie on the same geodesic from the base point of the distance function. This contradicts how geodesics must fill up balls under a $\MCP$ condition.

The strategy is similar to the proof of Theorem 4.3 of \cite{MMR}. Both arguments rely on the fact that the $\MCP$ condition forces geodesics to fill up small balls. However, we replace their use of \cite[Theorem 3.2]{HLD} on the tangent Carnot groups of $X$ with Cheeger differentiability of distance functions on $X$ itself. This is what allows us to drop the unique tangent hypothesis.

\subsection{Preliminaries}

Our space $(X,d,\m)$ is complete and separable, $\m$ is locally finite and Borel, and $\mathrm{supp} ~\m = X$.

A metric measure space $(X,d,\m)$ is {\em $n$-rectifiable} if $\m \ll \cH^n$ and there exists a countable family of Lipschitz functions $f_i$ defined on measurable subsets $A_i \subset \R^n$ so that $\m\left(X \setminus \bigcup_{i \in \N} f_i(A_i) \right) = 0$. If $n \notin \N$, then by convention, $X$ is not $n$-rectifiable.

A metric measure space $(X,d,\m)$ is a {\em Lipschitz differentiability space} if there exists a countable collection of Lipschitz maps $\{\phi_i : X \to \R^{n_i}\}_{i \in \N}$ and positive measure Borel sets $\{U_i\}_{i \in \N}$ so that $\m\left(X \setminus \bigcup_{i \in \N} U_i \right) = 0$ and for every Lipschitz $f : X \to \R$, $i \in \N$, and $\m$-a.e. $x \in U_i$, there exists a unique $Df(x) \in L(\R^{n_i}, \R)$ so that
\begin{align*}
  f(y) = f(x) + Df(x)(\phi_i(y) - \phi_i(x)) + o(d(x,y)).
\end{align*}
These were introduced (not with this name) by Cheeger in \cite{cheeger}. In the same paper, he proved that doubling metric measure spaces satisfying a Poincar\'e inequality (that we will not introduce here) are such spaces.

A geodesic in $X$ is a map $\gamma : [0,1] \to X$ for which
\begin{align*}
  d(\gamma(s),\gamma(t)) = |s-t| d(\gamma(0),\gamma(1)), \qquad \forall s,t \in [0,1].
\end{align*}
We let $\Geo(X)$ denote all the geodesics of $X$ and define $\ell(\gamma) := d(\gamma(0),\gamma(1))$. For $t \in [0,1]$, we let $e_t(\gamma) := \gamma(t)$ for $\gamma \in \Geo(X)$.
Let
\begin{align*}
  s_K(t) = \begin{cases}
    \tfrac{1}{\sqrt{K}} \sin(\sqrt{K}t), & K > 0, \\
    t & K = 0, \\
    \tfrac{1}{\sqrt{-K}} \sinh(\sqrt{-K}t), & K < 0.
  \end{cases}
\end{align*}
Let $v_{K,N}(r) = c_N\int_0^r s_K(t/\sqrt{N-1})^{N-1} ~dt$ and $\beta_t(r) = t\left(\frac{s_K(tr/\sqrt{N-1})}{s_K(r/\sqrt{N-1})}\right)^{N-1}$ where $c_N > 0$ is chosen so that $\lim_{r \to 0} \frac{v_{K,N}(r)}{r^N} = 1$. Note that $s_K(t) = t + O_K(t^3)$.

Following Ohta, for $N \in (1,\infty)$, we say a metric measure space $(X,d,\m)$ satisfies the $(K,N)$-measure contraction property or $\MCP(K,N)$ if for every $x \in X$ and measurable $A \subset X$ (or $A\subset B(x,\pi\sqrt{(N-1)/K})$ when $K > 0$) for which $0 < \m(A) < \infty$, there is a probability measure $\Pi$ on $\Geo(X)$ so that
\begin{enumerate}[label=(\arabic*)]
  \item We have $(e_0)_*\Pi = \delta_x$ and $(e_1)_*\Pi = \frac{1}{\m(A)}\m|_A$,
  \item For $t \in [0,1]$, $(e_t)_* \left( \beta_t(d(x,\gamma(1))) \m(A) ~d\Pi(\gamma) \right) \leq \m$.
\end{enumerate}

\subsection{AI use disclosure}
The author used Claude Fable 5.1 throughout the project. The intuition that lower dimensional Cheeger charts should contradict how geodesics fill up balls in $\MCP$ spaces had been known, but the set theoretic lim sup argument was proposed by Claude. Claude was also used to generate the figure, check for errors, and improve the presentation. The author wrote the paper, verified the math, and takes full responsibility for the contents of this paper.

\section{Proof of main theorem}
By Bishop-Gromov \cite[Theorem 5.1]{ohta}, $\frac{\m(B(x,r))}{v_{K,N}(r)}$ is non-increasing in $r$ for every $x$. A simple Taylor approximation calculation gives that $v_{K,N}(t) = t^N + O_{K,N}(t^{N+2})$. Thus, $\m(B(x,r)) \geq \frac{\m(B(x,R_0))}{v_{K,N}(R_0)}v_{K,N}(r) \gtrsim r^N$ for $r < R_0$ where
$$R_0 =
\begin{cases}
  \min\{1,\tfrac{1}{2}\pi\sqrt{(N-1)/K}\}, & K > 0, \\
  1, & K \leq 0,
\end{cases}$$
and so $\Theta^N_*(x) > 0$ for all $x$. We also get $\m(B(x,2r)) \leq \frac{v_{K,N}(2r)}{v_{K,N}(r)} \m(B(x,r))$ and as $\lim_{r \to 0} \frac{v_{K,N}(2r)}{v_{K,N}(r)} = 2^N$, $X$ is locally uniformly doubling.

By \cite[Theorem 1.3]{EB}, $X$ satisfies a Poincar\'e inequality and, together with locally uniformly doubling, is a Lipschitz differentiability space \cite{cheeger}. Thus, there is a collection of charts $\{\phi_i : X \to \R^{k_i}\}_{i \in \N}$ and positive measure Borel sets $\{U_i\}_{i\in\N}$ for which Lipschitz functions are $\phi_i$-differentiable for $\m$-a.e. $x \in U_i$ and $\m\left(X \setminus \bigcup_{i \in \N} U_i \right) = 0$. By the $N$-density bounds, we also get that $\m \ll \cH^N$.

Suppose $X$ is not $N$-rectifiable for contradiction. As $\Theta^N_*(x) > 0$ everywhere, $\dim_H (\m|_{U_i}) \leq N$ and so by \cite[Theorem 5.3]{BKO}, $k_i \leq N$ for all $i$.  We claim there is some $k_i < N$. If $N$ is non-integral, this is obvious. If $N$ is an integer, this follows from the density bounds, the fact that $X$ is not $N$-rectifiable, and \cite[Theorem 1.2]{BL}. We now choose such a $U_i$. Recalling that Bishop-Gromov gives positive lower $N$-density bounds, we let $F \subset U_i \cap \{\Theta^{N*} < \infty\}$ be a positive finite measure compact set so that there exists $0 < b < \infty$ and $R > 0$ so that
\begin{align}
  \m(B(x,r)) > br^N, \qquad \forall r < R, x \in F. \label{eq:AR}
\end{align}
Let $\phi = \phi_i : X \to \R^k$.

For any $p \in X$, $d_p(x) := d(x,p)$ is Lipschitz and so it is $\phi$-differentiable for $\m$-a.e. $x \in F$. By standard measure theoretic arguments $\{(x,p) \in F \times X: Dd_p(x) ~\mathrm{exists}\}$ is measurable and so by Fubini, we get that there is a full measure set $G \subset F$ so that for every $x \in G$, $d_p$ is differentiable at $x$ for $\m$-a.e. $p \in X$.

We claim for $\m$-a.e. $x \in G$, there is a sequence $\{x_j\}_{j\in \N} \subset G$ so that $x_j \to x$ and $|\phi(x) - \phi(x_j)| = o(d(x,x_j))$. Otherwise, there is some positive measure $G' \subset G$, $\kappa > 0$, and $r_0 > 0$ so that $|\phi(x) - \phi(y)| > \kappa d(x,y)$ for $x,y \in G'$ so that $d(x,y) < r_0$. But then $\phi|_{G' \cap B(x,r_0/2)}$ is biLipschitz into $\R^k$ so that $\cH^N(G' \cap B(x,R/2)) = 0$ for every $x \in G'$ as $k < N$. As $\m \ll \cH^N$, we get $\m(G') = 0$, a contradiction. Note then for $\m$-a.e. $y \in X$, we get
\begin{align}
  |d_y(x_j) - d_y(x)| \leq |Dd_y(x)(\phi(x_j) - \phi(x))| + o(d(x,x_j)) = o(d(x,x_j)). \label{eq:not-collinear}
\end{align}
Informally, $x_j$ and $x$ cannot be collinear with regards to geodesics from $y$.

Now choose $x \in G$ with such a sequence $x_j \to x$ and let $r_j = d(x,x_j)$. 
Let $\eta = 1/8$, $\theta := \lim_{r \to 0} \frac{\m(B(x,r))}{v_{K,N}(r)}$, $\delta > 0$ to be determined, and $\rho > 0$ also to be determined, but so that $\frac{\m(B(x,\rho))}{v_{K,N}(\rho)} > \theta - \delta$ and $\left| \frac{v_{K,N}(r)}{r^N} - 1 \right| \leq \frac{\delta}{\theta}$ for all $r \leq \rho$ (recall $0 < \theta < \infty$ by \eqref{eq:density-bounds} and \eqref{eq:AR}). Then
\begin{align}
  (\theta - 2\delta) r^N \leq (\theta - \delta) v_{K,N}(r) \leq \m(B(x,r)) \leq \theta v_{K,N}(r) \leq (\theta + \delta) r^N, \qquad \forall r \leq \rho. \label{eq:BG}
\end{align}

Note that $(1 - C\rho^2)^{N-1} \leq \frac{\beta_{t}(r)}{t^N} \leq (1+C\rho^2)^{N-1}$ for all $r < \rho$ and $t \in (0,1]$.
Let $j$ be large enough so that $t_j := (1+\eta) \frac{r_j}{\rho} < 1$. Let $\Pi$ be a probability measure on $\Geo(X)$ from $\delta_x$ to $\frac{1}{\m(B(x,\rho))}\m|_{B(x,\rho)}$ as dictated by $\MCP(K,N)$.
Define
$$\m_{t} = \m(B(x,\rho)) (e_{t})_* (\beta_t(d(x,\gamma(1))) d\Pi).$$
Then $\m_{t_j}$ is supported in $B_j := B(x,(1+\eta)r_j)$ and
\begin{align*}
  \m_{t_j}(B_j) &\geq \int \beta_{t_j}(d(x,\gamma(1))) \m(B(x,\rho)) ~d\Pi \\
  &\geq (1 - C\rho^2)^{N-1} t_j^N\m(B(x,\rho)) \\
  &\overset{\eqref{eq:BG}}{\geq} (\theta - 2\delta)(1-\tfrac{\delta}{\theta}) (1+\eta)^N r_j^N \\
  &\geq (\theta - 3\delta) (1+\eta)^N r_j^N
\end{align*}
when $\rho$ is sufficiently small.
Let $F_j := B(x_j,\eta r_j)$ and let $j$ be large enough so that $\eta r_j < R$. Then $F_j \subset B_j$ and so by $\MCP(K,N)$ (this is the only point $\m_t \leq \m$ is used),
\begin{align*}
  \m_{t_j}(B_j \setminus F_j) \leq \m(B_j \setminus F_j) \overset{\eqref{eq:AR} \wedge \eqref{eq:BG}}{\leq} (\theta + \delta)(1+\eta)^N r_j^N - b \eta^N r_j^N.
\end{align*}
Together, we get
\begin{align}
  \m_{t_j}(F_j) \geq (b\eta^N - 4\delta(1+\eta)^N) r_j^N \geq \frac{1}{2} b\eta^Nr_j^N \label{eq:big-F}
\end{align}
where the last inequality follows from choosing $\delta$ sufficiently small.

Disintegrating $\Pi$ along $e_1$ gives a measurable family of probability measures $\{\Pi_y\}_{y \in B(x,\rho)}$, each supported on geodesics between $x$ and $y$. For $0 < t < 1$, we define $\Phi_t(y,\cdot) := (e_t)_* \Pi_y$ so that,
$$\m_t(A) = \int_{B(x,\rho)} \beta_{t}(d(x,y)) \Phi_t(y,A) ~d\m(y).$$
Let $Z_j = \{y \in B(x,\rho): \Phi_{t_j}(y,F_j) > 0\}$. Then
\begin{align*}
  \m_{t_j}(F_j) &= \int_{B(x,\rho)} \beta_{t_j}(d(x,y)) \Phi_{t_j}(y,F_j) ~d\m(y) \\
  &\leq \int_{Z_j} \beta_{t_j}(d(x,y)) ~d\m(y) \leq t_j^N (1 + C\rho^2)^{N-1} \m(Z_j).
\end{align*}
Together with \eqref{eq:big-F}, we get
\begin{align*}
  \m(Z_j) \geq \frac{1}{2}b \eta^N r_j^N (1+C\rho^2)^{1-N}t_j^{-N} \gtrsim \rho^{N}.
\end{align*}

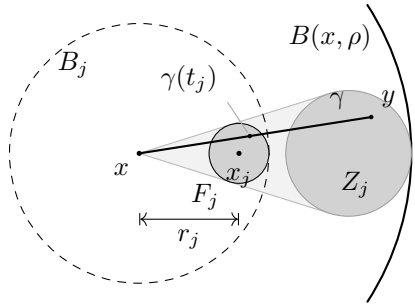
\begin{figure}[h!]
  \centering
\begin{tikzpicture}[scale=0.6, line cap=round, line join=round, font=\small]
  \def\rhoo{6}                             
  \def\rj{2.2}                             
  \def\eta{0.3}                            
  \pgfmathsetmacro{\Fr}{\eta*\rj}          
  \pgfmathsetmacro{\Bj}{(1+\eta)*\rj}      
  \pgfmathsetmacro{\tj}{\Bj/\rhoo}         
  \pgfmathsetmacro{\Zc}{\rj/\tj}           
  \pgfmathsetmacro{\Zr}{\Fr/\tj}
  \pgfmathsetmacro{\alp}{asin(\eta)}       
 
  \coordinate (x)  at (0,0);
  \coordinate (xj) at (\rj,0);
  \coordinate (Z)  at (\Zc,0);
 
  \fill[gray!10] (x) -- ($(Z)+({90+\alp}:\Zr)$)
      arc[start angle={90+\alp}, end angle={-90-\alp}, radius=\Zr] -- cycle;
  \draw[gray!55, thin] (x) -- ($(Z)+({90+\alp}:\Zr)$);
  \draw[gray!55, thin] (x) -- ($(Z)+({-90-\alp}:\Zr)$);
 
  \draw[thick] (-33:\rhoo) arc[start angle=-33, end angle=33, radius=\rhoo];
  \node[anchor=east] at (5.35,2.55) {$B(x,\rho)$};
 
  \draw[dashed] (x) circle (\Bj);
  \node at (-1.45,1.95) {$B_j$};
 
  \fill[gray!35] (Z) circle (\Zr);
  \draw[gray!75] (Z) circle (\Zr);
  \node at ($(Z)+(0.15,-0.7)$) {$Z_j$};
  \fill[gray!35] (xj) circle (\Fr);
  \draw (xj) circle (\Fr);
  \node at ($(xj)+(-128:\Fr+0.55)$) {$F_j$};
 
  \coordinate (y) at ($(Z)+(0.5,0.8)$);
  \coordinate (q) at ($\tj*(y)$);
  \draw[thick] (x) -- (y) node[pos=0.87, above=1pt, sloped] {$\gamma$};
  \fill (q) circle (1.4pt);
  \node (qlab) at (1.1,1.65) {$\gamma(t_j)$};
  \draw[gray!70, thin] (qlab.south east) -- ($(q)+(-0.05,0.08)$);
  \fill (y) circle (1.4pt) node[above right=0pt] {$y$};
 
  \fill (x) circle (1.6pt) node[below left=0pt] {$x$};
  \fill (xj) circle (1.6pt) node[below=2pt] {$x_j$};
  \draw[|<->|, thin] (0,-1.45) -- (\rj,-1.45) node[midway, below] {$r_j$};
\end{tikzpicture}
  \caption{How $F_j$ and $Z_j$ are positioned relative to $x$ in $B(x,\rho)$}
\end{figure}

For $\m$-a.e. $y \in Z_j$, there exists a geodesic $\gamma$ from $x$ to $y$ so that $q = \gamma(t_j) \in F_j$. Thus, $d(x,q) + d(q,y) = d(x,y)$ and so
\begin{align}
  d_y(x) - d_y(x_j)  = d(x,y) - d(x_j,y) \geq d(x,q) - d(x_j,q) \geq (1 - 2\eta)r_j. \label{eq:collinear}
\end{align}

Finally, let $Z = \limsup_j Z_j$. Since $Z_j \subset B(x,\rho)$ and $\m(B(x,\rho)) < \infty$, we get that $\m(Z) \geq \limsup_j \m(Z_j) \gtrsim \rho^{N} > 0$. Then for $\m$-a.e. $y \in Z$, \eqref{eq:collinear} is satisfied for an infinite subsequence $x_{j_\ell} \to x$. But this contradicts \eqref{eq:not-collinear}.

\end{document}